\documentclass{article}
\usepackage[utf8]{inputenc}
\usepackage{amsmath}
\usepackage{amssymb}
\usepackage{graphicx}
\usepackage{amsthm}
\usepackage{float}
\usepackage{fancyref}

\graphicspath{ {./images/} }

\title{Subsums of Null Sequences: A Survey and Open Questions }
\author{Justin Jacob }
\date{\today}

\begin{document}
\maketitle
\begin{center}
  \textbf{Abstract}
\end{center}
We give a survey on the different results involving the topological structure of subsums of null sequences. 
\section{Achievement sets and Null Sequences:}
Consider a convergent sequence $x(n)=[x_0,x_1,{\ldots},x_n,{\ldots}],$ such that $\displaystyle \lim_{n \to \infty}{x_n}=0$. Denote the set of subsums of this convergent sequence as: $$E(x)=\sum_{i=1}^{\infty} \epsilon_{i}x_i, \epsilon_{i} \in \left\{0,1\right\}, \forall i \in \mathbb{N}$$
where $E(x)$ called the $\textit{achievement set}$ of $x(n)$. Define the $n$-th tail as: $$R_{n} =\sum_{i=n+1}^{\infty} x_i$$
For example, we can see that the sequence $x(n)=[ \frac{1}{2}, \frac{1}{4},{\ldots}, \frac{1}{2^{n}},{\ldots}]$ is a convergent sequence. Notice that $\sum_{n=1}^{\infty} x(n) = 1$. Thus we know the set $E(x) \subseteq [0,1]$. To show $E(x)$ is equal to $[0,1]$ we use the following proof from Nitecki \cite{nitecki2013subsum}. Observe that to express a digit in base-$2$, or binary, we can use the following notation. Consider $\epsilon_{i} \in \{0,1\}, \forall i \in \mathbb{N}$. Then for any $x \in E(x)$ we have that $x= \sum_{i=1}^{\infty} \frac{\epsilon_{i}}{2^i}$ for some sequence $[\epsilon_0,\epsilon_1,\epsilon_2,{\ldots}]$ Since every number between 0 and 1 has a binary expansion we can see that $E(x)=[0,1]$.
\newline
For another example, consider the sequence $x(n)=[ \frac{1}{3}, \frac{1}{9},{\ldots}, \frac{1}{3^{n}},{\ldots}].$ Again take $\epsilon_{i} \in \{0,1\}, \forall i \in \mathbb{N}$. We first give the definition of the standard Cantor set.
\newtheorem{definition}{Definition}
\begin{definition}[Standard Cantor Set]
The standard Cantor set $C$ is defined as removing the middle third from the interval $[0, 1]$, removing the middle third from each of the remaining connected components,  and so on, iteratively continuing this process. Thus: $$C=[0,1] \setminus \bigcup\limits^{\infty}_{k=1}{C_k},  \text{where} \; C_k=\frac{C_{k-1}}{3}\cup \left(\frac{2}{3}+\frac{C_{k-1}}{3}\right), n \geq 2  \; \text{and}\; C_1=\left[\frac{1}{3}, \frac{2}{3}\right].$$
\end{definition}
\newtheorem{theorem}{Theorem}
\begin{theorem}[Nitecki \cite{nitecki2013subsum}]
The set $E(x)$ for the sequence $x(n) = \frac{1}{3^n}, n \in \mathbb{N}$ is homeomorphic to the standard Cantor set.
\begin{proof}
First notice that any subsum must belong to the interval $[0,1/2]$. However, unlike the previous example, every sum excluding the first term will belong in the interval $J_{0}=[0,\frac{1}{6}]$, while $J_{1}=[\frac{1}{3},\frac{1}{2}].$ So thus subsums are contained in $C_1=J_{0} \cup J_{1}$. Now we consider the exclusion of the second term, $\frac{1}{9}$, we find that the subsum is contained in four subintervals, two subintervals of $J_0$,$J_{00}$ and $J_{01}$ and two subintervals of $J_1$, $J_{10}$ and $J_{11}$. Notice these correspond to all possible sequences of $\epsilon_i$ of length 2. Then continuing this process for all possible sequences of $\epsilon_i$ for length $n$, or where $1 \leq i \leq n$, we have that:
$$E(x) \subset C_n = {\bigcup_{{\epsilon_{1},\epsilon_{2},{\ldots},\epsilon_{n} \in \{0,1\}^{n}}}^{}} J_{\epsilon_{1},{\ldots},\epsilon_{n}}$$ 
With $J_{0^n}=[0,R_{n}]$. Notice that at the next step, we have two possibilities for the next subintervals: $${J_{{\epsilon_{1},{\ldots},\epsilon_{n},0}}} \,\text{and} \,  {{J_{\epsilon_{1},{\ldots},\epsilon_{n},0}}}+\textstyle\frac{1}{3^{n+1}}= {{J_{\epsilon_{1},{\ldots},\epsilon_{n},1}}}$$
So at every step, our intervals are disjoint. We can notice that passing from the union $c_n$ to $C_{n+1}$, a gap is created. Furthermore notice that $R_{n+1}=\frac{1}{3} R_n$, so at each step the size of this gap is the middle third of each component. Thus we have $C_{\infty}= \bigcap_{n=1}^{}C_n$, which is a version of the standard Cantor set, just scaled down by a factor of $\frac{1}{2}$. In our construction, our Cantor set contains every number whose ternary expansion contains digits of only zeroes and ones rather than the zeroes and twos in the standard Cantor set.
\end{proof}
\end{theorem}

\begin{figure}[h!]
 \caption{The construction of subsum sets for  $x(n) = \frac{1}{3^n}, n \in \mathbb{N}$ \cite{aguero_2019}}

\includegraphics{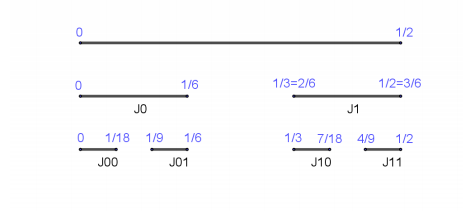}
\end{figure}
So far for the current examples considered, one may notice that in our investigation of the achievement sets for the geometric series $x(n) = \frac{1}{3^n}$ and $x(n) = \frac{1}{2^n}, n \in \mathbb{N}$, either the intervals at each step will continue to overlap or be completely disjoint.
\newline
\par We can now generalize to obtain the following results about subsums of geometric series. The proof of the theorem found in Aguero's paper \cite{aguero_2019} is very similar to Nitceki's proof. The key insight is that if $r \geq \frac{1}{2}$, we will have that the intervals constructed at each step will never be completely disjoint and thus when we take the intersection we will have that $C_{\infty}= \bigcap_{n=1}^{}C_n$ must be an interval. However if $r <\frac{1}{2}$ then each interval will indeed be disjoint and thus $C_{\infty}= \bigcap_{n=1}^{}C_n$ will be homeomorphic to a Cantor set. Of course for a more rigorous notion one can consult \cite{aguero_2019}.
\newline
\begin{theorem}[\cite{aguero_2019}]
For $0 \leq r < \frac{1}{2}$ we have $E(x)$ for $x_k = r^{k}$ is a Cantor set and $ \frac{1}{2} \leq r <1$ $E(x)$ is an interval.
\end{theorem}

\begin{figure}[!htb]
 \caption{Subsum sets for $ \frac{1}{2} \leq r <1$ \cite{aguero_2019}}
\centering
\includegraphics{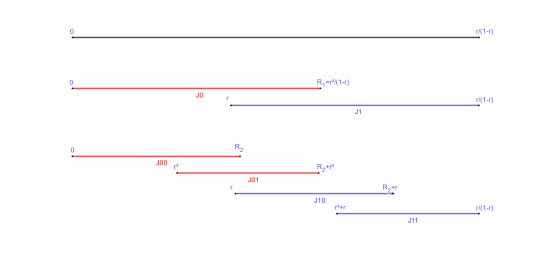}
\end{figure}
\section{General Results on $E(x)$ and Cantorvals}
We now shift our attention towards more general results regarding $E(x)$ for arbitrary null sequences, not just geometric ones. The first prominent result was done in 1914 by Kakeya \cite{Kakeya} showed that $E(x)$ is:
\begin{enumerate}
    \item A perfect set,
    \item The union of finite intervals if $x_n \leq R_n$ for $n$ sufficiently large.
    \item Homeomorphic to a Cantor set if $x_n >R_n$ for $n$ sufficiently large.
\end{enumerate}
Kakeya also conjectured that $E(x)$ is nowhere dense and homeomorphic to the Cantor set if $x_n > R_n$ for infinitely many $n$. However this conjecture has been shown to be false with many such counterexamples. Ferens \cite{Ferens} was the first to give the following example with proof:
$$\sum_{i=1}^{\infty} \left[\epsilon_{5i-4} \left(\frac{7}{27^i} \right) + \epsilon_{5i-3} \left(\frac{6}{27^i}\right) + \epsilon_{5i-2} \left(\frac{5}{27^i}\right) + \epsilon_{5i-1} \left(\frac{4}{27^i}\right) + \epsilon_{5i} \left(\frac{3}{27^i}\right)\right]$$
with $\epsilon_i \in \left\{0,1\right\}, \forall i \in \mathbb{N}$ as usual. Gutherie and Nymann \cite{Nymann} gave a much simpler example in 1988:
\begin{equation} \label{eq:solve} 
T=\left\{\sum_{i=1}^{\infty} \left[\epsilon_{2i-1} \left(\frac{3}{4^i} \right) + \epsilon_{2i} \left(\frac{2}{4^i}\right) \right]\right\}    
\end{equation}

\begin{figure}[h!]
 \caption{The set of subsums of $T=\left\{\sum_{i=1}^{\infty} \left[\epsilon_{2i-1} \left(\frac{3}{4^i} \right) + \epsilon_{2i} \left(\frac{2}{4^i}\right) \right]\right\}$}
\centering
\includegraphics{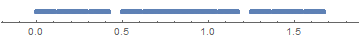}
\end{figure}
Gutherie and Nymann were also able to come up with a very important theorem regarding the classification of sets of $E(x)$ for arbitary null sequences. Before revealing their result, we define a Cantorval, dubbed by Mendes and Oliviera in the following paper \cite{Mendes}. However, we give more comprehensive definition regarding them than the one in their paper:

\begin{definition}[Cantorval \cite{anisca_chlebovec_ilie_2012}]
Let $K$ be a compact subset of $\mathbb{R}$. A gap of $K$ is a bounded connected component of $\mathbb{R} \setminus K$; an interval of $K$ is a non-trivial connected
component of $K$. A perfect subset of $\mathbb{R}$,  such that any gap is accumulated on each side by infinitely many intervals and gaps is called an M-Cantorval. A perfect subset of $\mathbb{R}$,  such that any gap has an interval adjacent on its right and is accumulated on the left by infinitely many intervals and gaps is called an L-Cantorval. The definition of an R-Cantorval is analogous.
\end{definition}
An analgous definition of an M-Cantorval is also any subset of $\mathbb{R}$ that is homeomorphic to $C \cup \left( \bigcup_{n=1}^{\infty} C_{2n-1} \right)$.
Observe that when dealing with Cantorvals we have the following set equality: 
$$C \cup \left( \bigcup_{n-1}^{\infty} C_{2n-1} \right) = C_{2n} \setminus \left( \bigcup_{n-1}^{\infty} C_{2n} \right)$$
Many of the topological properties of Cantorvals have been considered in the literature, see \cite{banakh2012algebraic} for some of these results. For special types of M-Cantorvals found in \cite{banakiewicz2017},their Lebesgue measures have been shown to be the sum of the measures of their interval components. However, Cantorvals do not only occur in the topic of achievement sets. Bartoszewics et. al. \cite{Bartoszewicz_2018} showed that Cantorvals can occur as attractors of iterated function systems consisting of affine functions.
\newline
We state a lemma attributed to Nymann and Saenz \cite{Nymann}:
\newtheorem{lem}{Lemma}
\begin{lem}
If (a,b) is a gap of $E(x)$, then for some $\epsilon >0$ and $\epsilon' >0$,
$$b+([0,\epsilon] \cap E) =[e,b+\epsilon] \cap E \; \text{and} \; [s-\epsilon',s] \cap E = (s-a) +([a-\epsilon',a] \cap E)$$
\end{lem}
\par
Now we are ready to state the Gutherie-Nymann-Saenz Classification Theorem:
\begin{theorem}[Gutherie-Nymann Saenz \cite{Nymann}]
The achievement set $E(x)$ for a general null sequence $x(n)=[x_0,x_1,{\ldots},x_n,{\ldots}]$ is homeomorphic to one of the following:
\begin{enumerate}
    \item A Cantor set,
    \item A union of finite intervals,
    \item Homeomorphic to the set $T$ defined in equation~\ref{eq:solve} (which is also an M-Cantorval).
\end{enumerate}
\begin{proof}
\par We have seen examples of null sequences whose achievement sets $E(x)$ are either a union of finite intervals or a Cantor set. Now suppose that E is neither a finite union of intervals nor homeomorphic to the Cantor set. Then it is clear the complement of E must contain infinitely many intervals. E must contain infinitely many intervals as well, for there were only finitely many then either $E \cap \left[0,\epsilon\right)$ is an interval or $E \cap \left[0,\epsilon\right)$ contains no interval for some $\epsilon >0$. If the former is true, then $\exists n$ such that $R_n$ is an interval and thus by Kakeya's result $E$ is a finite union of intervals. In the latter case, $\exists n$ such that $R_n$ is a Cantor set, so $E$ would be a Cantor set. In either case we arrive at a contradiction, so $E$ must have infinitely many intervals. In fact $E \cap [a,b]$ cannot be homeomorphic to the Cantor set $\forall a,b \in E$ since every tail $R_n$ has intervals in its set of subsums. Then for $x \in E$ we have $E \cap (x,x+\epsilon) \neq \emptyset$, and since $E$ is perfect we have $E \cap (x-\epsilon,x) \neq \emptyset$ for some $\epsilon >0$. Thus we can see that there are intervals dense in $E$.
\par Now we define a strictly increasing map $f$ from the union of all intervals in $T$ to the union of all intervals in $E$. We do this inductively: first map the longest interval in $T$ in a strictly increasing way to the longest intervals in $E$. There can be at most finitely many intervals of the same length in either set, so if no one interval is the longest we proceed by choosing the leftmost one. 
\par After the $n$-th step, $\frac{4^{n}-1}{3}$ intervals of $T$ will have been identified, in a strictly increasing way, with $\frac{4^{n}-1}{3}$ intervals of $E$. In addition, $\frac{2(4^{n}-1)}{3}$ gaps of $T$ will be identified in a strictly increasing way onto $\frac{2(4^{n}-1)}{3}$ gaps of $E$. 
We repeat the inductive process above to each of the spaces between any two adjacent intervals and gaps (of which there are $4^{n}-1$) and to the space between 0 and the leftmost of the gaps and to the space between the rightmost of the gaps and $\frac{5}{3}$. Lemma 1 guarantees this is possible, since there will be infinitely many intervals and gaps in $T$ and $E$ in $[0,\epsilon]$ and $[\frac{5}{3}-\epsilon,\frac{5}{3}]$, respectively.
\par
When $f$ is defined this way, it is a strictly increasing mapping from the union of all intervals in $T$ to the union of all intervals in $E$. As $T$ and $E$ are both perfect it is possible to extend $f$ continuously to all of $T$ and then onto all of $E$. Since the extension is strictly increasing, and $T$ is compact, $f$ is the desired homeomorphism.
\end{proof}
\end{theorem}
\begin{theorem}[\cite{10.4169/amer.math.monthly.122.9.862}]
Any two M-Cantorvals (symmetric Cantorvals) are homeomorphic.
\end{theorem}
We introduce the following notation introduced by Bartoszewicz in \cite{bartoszewicz2013multigeometric}:

\begin{definition}\label{def:notation}
For any $ 0 \leq q < \frac{1}{2}$ we use $(k_1,k_2,{\ldots},k_m;q)$ to denote the sequence:
\newline
$$\left(k_1,k_2,{\ldots},k_m,k_{1}q,k_2q,{\ldots},k_{m}q,k_1{q^2},k_2{q^2},{\ldots},k_m{q^2},{\ldots}\right)$$
\end{definition}

In 2011, Jones \cite{Jones2011AchievementSO} gave a result that yielded a continuum of Cantorvals:
\begin{theorem}
The set of subsums of $\left(\frac{3}{5},\frac{2}{5},\frac{2}{5},\frac{2}{5};q\right)$
is a Cantorval for all $q$ such that $$ \frac{1}{5} \leq \sum_{i=1}^{\infty} q^i < \frac{2}{9}$$
\end{theorem}
Kenyon and Nitecki \cite{10.4169/amer.math.monthly.122.9.862} also paved the way for another continuum of Cantorvals with the following result in 2015:
\begin{theorem}
Suppose we are given $n \in \mathbb{N}$ and $n$ integers $d_0,d_1, {\ldots},d_{n-1}$ such that $d_i \equiv i \pmod n$.
Then the set of generalized base $n$ expansions using these ``digits'':
$$ S=\left\{\sum_{i=1}^{\infty} \frac{a_i}{n^i} \, | \, {a_i} \in \left\{d_0,d_1, {\ldots},d_{n-1}\right\} \right\}$$
has nonempty interior.
\end{theorem}

Using the notation described in Definition \ref{def:notation}, Bartoszewics \cite{bartoszewicz2013multigeometric} established the following result:
\begin{theorem}
Let $k_1 \leq k_2 \leq \cdots \leq k_m$ be positive integers and $K=\sum_{i=1}^{m}{k_i}$. Assume that there exist positive integers $n_0$ and $n$ such that each of the numbers $n_0,n_{0}+1, {\ldots},n_{0}+n$ can be obtained by summing the numbers $k_1,k_2, {\ldots},k_m$ in some combination. Then we have:
\begin{enumerate}
    \item If $q \geq \frac{1}{n+1}$ then $E(k_1,k_2,{\ldots},k_m;q)$ has nonempty interior.
    \item If $q < \frac{k_m}{K+{k_m}}$ then $E(k_1,k_2,{\ldots},k_m;q)$ is not a finite union of intervals.
    \item If $\frac{1}{n+1} \leq q <\frac{k_m}{K+{k_m}}$ then $E(k_1,k_2,{\ldots},k_m;q)$ is a Cantorval.
\end{enumerate}
\end{theorem}
Ferdinands et. al. \cite{familyCantorvals} have generalized this result to the following:
\begin{theorem}
Let $k_1 \leq k_2 \leq \cdots \leq k_m$ be positive integers and $K=\sum_{i=1}^{\infty}{k_i}$. Assume that there exist positive integers $a$,$d$ and $n$ such that the set: $$\sum_{i=1}^{m} \epsilon_{i}k_i, \epsilon_{i} \in \{0,1\}$$ contains each of the numbers $a,a+d,a+2d,{\ldots},a+nd$. Then we have:
\begin{enumerate}
    \item If $q \geq \frac{1}{n+1}$ then $E(k_1,k_2,{\ldots},k_m;q)$ has nonempty interior.
    \item If $q < \frac{k_m}{K+{k_m}}$ then $E(k_1,k_2,{\ldots},k_m;q)$ is not a finite union of intervals.
    \item If $\frac{1}{n+1} \leq q <\frac{k_m}{K+{k_m}}$ then $E(k_1,k_2,{\ldots},k_m;q)$ is a Cantorval.
\end{enumerate}
\end{theorem}
Notice, using Definition \ref{def:notation}, we have the Gutherie-Nymann example $T=(3,2,\frac{1}{4})$.
Agüero \cite{aguero_2019} also generalized the Gutherie-Nymann example as following:
\begin{theorem} 
Using the notation in Definition~\ref{def:notation}, 
$(\underbrace{{3,3,{\ldots}}}_{k \; 3's},\underbrace{{2,2,{\ldots}}}_{k \; 2's},\frac{1}{3k+2k'+1})$ is a Cantorval for $k,k' \in \mathbb{N}$.
\end{theorem}
\section{Sums of Achievement Sets:}
We now turn our attention towards sums of achievement sets, ie: $$P=E_{m}(x)=\sum_{i=1}^{\infty} \epsilon_{i}x_i, \epsilon_{i} \in \{0,1,2,{\ldots},m\}, \forall i \in \mathbb{N}$$ with tail: $$R_z =\sum_{i=z+1}^{\infty} x_i $$
In 1995, Nymann \cite{Nymann1995} gave a result regarding these sums:
\begin{theorem}
There is a positive integer $m$ for which $P$ is a finite union of intervals if and only if $\limsup{\frac{x_n}{R_n}}<\infty$.
\newline
Moreover, the smallest positive integer for which $P$ is a finite union of intervals is the smallest integer $m$ such that $\frac{x_n}{R_n} \leq m$ for all but a finite number of integers $n$.
\end{theorem}
We now give a result that uses the proof Gutherie and Nymann used in their classification theorem:
\begin{theorem}[Generalized Guthrie-Nymann]
\label{Generalized Guthrie-Nymann}
 For all $ m \in \mathbb{N}$, $E_m(x)$ is homeomorphic to one of the following:
\begin{enumerate}
    \item A Cantor set,
    \item A union of finite intervals,
    \item A Cantorval.
\end{enumerate}
\begin{proof}
First notice that $E_m(x)$ is a perfect compact set since it is a sum of perfect compact sets. Then if the interior of ${E_m(x)} = \emptyset$, $E_m(x)$ must be homeomorphic to the Cantor set. 
\newline
Now suppose the interior of $E_m(x)$ is nonempty. Then all but finitely many tails $R_{n}$ must contain an interval as its set of subsums. Furthermore, we also have that $E(x) \subseteq E_m(x)$, so in this case we know at the very least $E_m(x)$ will contain a finite union of intervals. Now, we can use the same proof given in \cite{Nymann} to conclude the result.
\end{proof}
\end{theorem}
\section{Conjectures and Open Questions}
\subsection{Achievement Sets of Compact Sets with Decreasing Diameter:}
Now we consider a much more general problem: Define:
\newline
$${diam(A)}= \displaystyle \sup_{x,y \in A}\{d(x,y)\}$$ where $d(x,y)$ denotes the usual metric on $\mathbb{R}$. Suppose we are given a sequence of compact sets in $\mathbb{R}$, $K_1,K_2,{\ldots},K_i,{\ldots}$ where $0 \in K_i \, \forall i$, $\forall x \in K_i \, x \geq 0$ and ${diam(K_i)}>0 \; \forall i \in \mathbb{N}$. Furthermore let ${diam(K_i)} \to 0$ as $i \to \infty$. Now pick $x_1 \in K_1, x_2 \in K_2, {\ldots},x_n \in K_n, {\ldots}$ and consider: $$S=\{\sum_{j=1}^{\infty}x_j| \,x_j \in K_j\} \, \, \text{with tail} \, \, R_z=\{\sum_{i=z+1}^{\infty} x_i| \, x_i \in K_i\}$$

\newtheorem{lemma}{Lemma}
\begin{lemma}
S is a perfect set.
\begin{proof}
 Clearly by Tychnoff's Theorem $S$ is compact. Pick $x \in S$. Then we have $x=x_1+x_2+{\ldots}+x_n+{\ldots}$. Notice that since $diam(K_i) \to 0$ as $i \to \infty$, $\forall \epsilon >0, \, \exists K_y $ with $diam(K_y) <\epsilon$. Then $z_y \in K_y$ has the property that $\left\vert z_{y}-x_{y} \right \vert < \epsilon$, so the sum $x_1+x_2+{\ldots}+x_n+{\ldots}+z_y+{\ldots} \in (x,x+\epsilon)$. Then pick $z_z \in K_z$ where $diam(K_z)<\epsilon_{1}<\epsilon$, so $\left\vert z_{z}-x_{z}\right\vert < \epsilon_{1}$. So then $x_1+x_2+{\ldots}+x_n+{\ldots}+z_y+z_z+{\ldots} \in (x,x+\epsilon_{1}).$ Inductively doing this approach, we get a sequence in $S$ with limit point $x$. Thus $S$ is perfect.
\end{proof}
\end{lemma}
However, it is unknown what the structure of $S$ is in general. It is trivial to show that if the interior of $S$ is empty, then $S$ must be a Cantor set, however we believe that the classifying these sets completely in terms of the tail $R_z$ is far too general of a problem. Instead, we list the following open questions:
\begin{enumerate}
    \item What are the possible sets $S$ can be homeomorphic to?
    \item Under what conditions on the tail $R_z$ will this set be an interval or a finite union of intervals?
    \item Under what conditions on the tail $R_z$ will this set be an M-Cantorval? An L or R Cantorval?
    
\end{enumerate}
\subsection{Achievement Sets of Cantor Sets:}
\par Sums of finite Cantor sets appear frequently in the literature and show up in spectral theory and homoclinic bifuractions. In addition, many of their topological properties have been considered for sums of more general types of  Cantor sets. Thanks to the Newhouse Gap Lemma \cite{newhouse_1974},  we also know if $\tau(C_{\alpha})\tau(C_{\beta})>1$,  then the sum of the two Cantor sets contain an interval. Here $\tau(C)$ denotes the thickness of the set,  see \cite{newhouse_1974} for a more precise notion. It is well known (see chapter 4 in \cite{kraft_1996}) that if:
\begin{equation} \label{eq:cantorset}
  \frac{\log{2}}{\log{(\frac{1}{\alpha})}}+\frac{\log{2}}{\log{(\frac{1}{\beta})}}<1 \Rightarrow \dim_{H}(C_{\alpha})+\dim_{H}(C_{\beta})<1
\end{equation}
then the Lebesgue measure of the set $C_{\alpha}+C_{\beta}$ is zero and  the sum be a Cantor set. Here $\dim_{H}$ denotes the Hausdorff dimension of a set. By result \ref{eq:cantorset} and the Newhouse Gap Lemma \cite{newhouse_1974},  we have knowledge of many of the sums of Cantor sets depending on the parameters of $\alpha, \beta.$ However, many open questions still remain. These results reveal a ``mysterious region" in the space of parameters $\alpha, \beta.$ We denote the area where the sum is an interval as ``large", the area where the sum has Lebesgue measure zero as ``small", and the ``mysterious region" as R in the figure below. Solomyak \cite{solomyak_1997} found that if:
 $$\dim_{H}(C_{\alpha})+\dim_{H}(C_{\beta})>1, $$the set $C_{\alpha}+C_{\beta}$ has positive Lebesgue measure for almost everywhere $\alpha, \beta \in (0, \frac{1}{2}).$ Furthermore, he and Peres also found this is true for any compact set $K$ and any middle $\alpha$ Cantor set $C_{\alpha}$ \cite{10.2307/117688}. This implies that for almost every parameter in the ``mysterious region,"  the sum of the corresponding homogeneous Cantor sets has positive Lebesgue measure.
\begin{figure}[h!]
 \caption{A graph of the space of parameters and the "Mysterious Region" R. \cite{solomyak_1997}}
\centering
\includegraphics{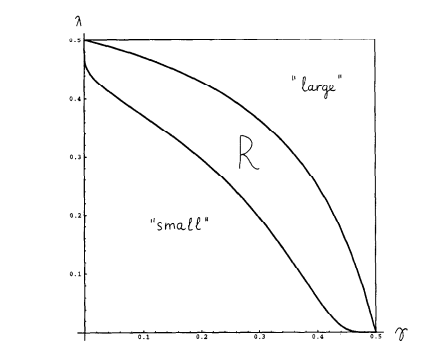}
\end{figure} 
 \par All of these results,however, concern finite sums of Cantor sets. In addition, many of the results are quite specific and cannot be verified in full generality, i.e. for every possible class of Cantor sets. However we believe many of the results will still hold even in our generalization of the problem: Suppose we have a null sequence $\alpha_1,\alpha_2,{\ldots},\alpha_n,{\ldots}$ such that $\alpha_i \to 0$ as $i \to \infty$. Furthermore suppose we have a sequence $\lambda_1,\lambda_2,{\ldots},\lambda_n,{\ldots}$ with $\lambda_i \to 0$ as $i \to \infty$. For all $i$ let $K_i$ be a middle-$\alpha_i$ Cantor set with first step $[0,\lambda_i]$. Now consider $S$ for these sets. 
We pose the following open questions when considering $S$:
\begin{enumerate}
    \item Under what conditions of $R_z$ and $\alpha_i$ for $i \in \mathbb{N}$ will this set be a Cantor set? What type of Cantor set will it be?
    \item Under what conditions of $R_z$ and $\alpha_i$ for $i \in \mathbb{N}$ will this set be a finite intervals or a union of finite intervals?
    \item Under what conditions of $R_z$ and $\alpha_i$ for $i \in \mathbb{N}$ will this set be an M-Cantorval? An L or R Cantorval?
    \item Denote the Haussdorf dimension of $S$ as $\dim_{H}(S)$. Is it true that $\dim_{H}(S) =\sum_{i=1}^{\infty} \dim_{H}(K_i)$ generically? If not, in what specific cases will this be true?
\end{enumerate}
\subsection{Achievement Sets in the Complex Plane:}
\par Our last generalization of the problem involves making each set $K_i$ a subset of the complex plane and then investigating it. To do this we give a brief overview of previous results relating to the problem at hand.
\par In 1957, Rényi \cite{renyi_1957} was the first to consider representations of numbers with radix in $\mathbb{R}$, i.e. the set of points that can be represented as $$x=\sum_{i=1}^{\infty} \frac{\epsilon_i}{q^i}$$ with $q \in \mathbb{R}$ and $\epsilon_k \in A$, where A is some finite set of positive real values called the $\textit{alphabet}$. Call this set $H$. This led to many more developments, including the study of canonical number systems and radix representations in algebraic number fields. See \cite{akiyama_borbely_brunotte_petho_thuswaldner_2004},\cite{akiyama_thuswaldner_2006} and \cite{akiyama_2005} for a more comprehensive summary. The first people to study complex radixes were Daróczy and Kátai \cite{dk88}, who showed for which $q^i=\theta^i$ for some $\theta \in \mathbb{C}$ with $\left\vert \theta \right\vert <1$ does $H= \mathbb{C}$. These conditions were expanded upon by Komornik and Loreti in \cite{komornik_loreti_2007}. The study of these complex bases has found various applications in computer science and cryptography, see \cite{gilbert_1984} and \cite{solinas_2000}. Gilbert \cite{gilbert_1981} was the first to consider the fascinating geometry of these expansions in 1981. One such example of these fractals is the Lévy dragon. See \cite{levydrag} for a very thorough analysis on this figure.
\begin{figure}[h!]
 \caption{The Lévy Dragon \cite{levydrag}}
\centering
\includegraphics{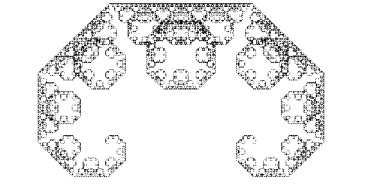}
\end{figure} 
Lai \cite{article} studied the geometry of these sets with base $q=pe^{\frac{2\pi i}{n}}$ and characterized their convex hull. Furthermore, the topological properties of these structures, called self-affine tiles, were done by Akiyama and Thuswaldner \cite{akiyama_thuswaldner_2004}. There are many open questions regarding this subject already posed in the literature, mostly regarding the representability of numbers in this radix system as well as the properties of the sets they form.
\par Finally we give the generalization of our problem: Again, define:
$${diam(A)}= \displaystyle \sup_{x,y \in A}\{d(x,y)\}$$ where $d(x,y)$ denotes the usual metric on $\mathbb{C}$. Suppose we are given a sequence of compact sets in $\mathbb{C}$, $K_1,K_2,{\ldots},K_i,{\ldots}$ where $0 \in K_i \, \forall i$, $\forall x \in K_i \, x \geq 0$,  ${diam(K_i)}>0 \; \forall i \in \mathbb{N}$ and ${diam(K_i)} \to 0$ as $i \to \infty$. Now pick $x_1 \in K_1, x_2 \in K_2, {\ldots},x_n \in K_n, {\ldots}$ and consider: $$S=\left\{\sum_{j=1}^{\infty}x_j| \,x_j \in K_j\right\} \, \, \text{with tail} \, \, R_z=\left\{\sum_{i=z+1}^{\infty} x_i| \, x_i \in K_i\right\}$$
Now we pose the following questions, inspired from the results regarding complex radix expansions and the sets they form:
\begin{enumerate}
    \item Under what conditions on $K_i$ and $R_z$ will $S=\mathbb{C}$?
    \item Several authors have investigated the fractal structures of radix representations using a complex base \cite{gilbertfractal}. Under what conditions on $K_i$ and $R_z$ will $S$ have fractal structure? What will the Haussdorf dimension of $S$ be? 
    \item Motivated by results in algebra concerning representations in algebraic number fields \cite{akiyama_borbely_brunotte_petho_thuswaldner_2004}, besides the trivial case, what specific conditions on $K_i$ and $R_z$ will we able to have a radix representation of all algebraic numbers?  
    \item In the case of real numbers, we can see Cantor sets and Cantorvals arise as achievements sets of null sequences. Besides the trivial case of real numbers already examined, under what conditions of $R_z$ and $\alpha_i$ for $i \in \mathbb{N}$ will $S$ be a Cantor subset of $\mathbb{C}$? A 2-dimensional Cantor set (Cantor Dust)?
    \item Motivated by the results in \cite{akiyama_thuswaldner_2004}, under what conditions of $R_z$ and $\alpha_i$ for $i \in \mathbb{N}$ will this set be a disk or a union of finite disks?
    \item Denote the Haussdorf dimension of $S$ as $\dim_{H}(S)$. Is it true that $\dim_{H}(S) =\sum_{i=1}^{\infty} \dim_{H}(K_i)$ generically? If not, in what specific cases will this be true?
\end{enumerate}
\section{Acknowledgements}
This project was done for the Summer Undergraduate Research Experience (SURP) Program at UC Irvine. I am grateful to the Undergraduate Research Opportunities Program (UROP) at UCI for providing a stipend to support the project. In addition, I would like to thank my advisor, Professor Anton Gorodetski, for guiding me throughout the entire process and supporting me. I would also like to thank Lucas Booher for his contributions towards the project as well.
\bibliography{message.bib} 
\bibliographystyle{plain}

\end{document}